\documentclass[proc]{edpsmath}

\usepackage[T1]{fontenc}
\usepackage{bbm}
\usepackage{graphicx}

\usepackage{afterpage}
\usepackage{url}

\usepackage{amssymb}
\usepackage{mathrsfs}
\usepackage{stmaryrd}

\newcommand	{\indic}	[1] {{\mathbbm{1}}_{#1}}
\newcommand	{\indicbis}	[1] {\indic{\{#1\}}}

\newcommand	{\dtf}	{\bar{d}_f}
\newcommand	{\dtfp}	{\bar{d}_{F_p}}

\newcommand{\tr}{\mathbbm{t}}
\newcommand{\trT}{\tr^{\{T\}}}

\providecommand{\Fonction}[5]
{\begin{alignat*}{2}
  & #1 :\; & #2 & \longrightarrow #3 \\
  & & #4 & \longmapsto #5
\end{alignat*}}

\newcommand{\E}{\mathbbm{E}}

\newcommand{\R}{\xR}

\newcommand{\N}{\xN}

\newcommand{\mc}{\mathcal{M}}

\begin{document}

\title{Random ultrametric trees and applications} 
\author{Amaury Lambert}\address{Universit\'e Pierre et Marie Curie -- UPMC Univ Paris 06,
 Laboratoire de Probabilit\'es et Mod\`eles Al\'eatoires,
 Case 188, 
 4 place Jussieu,
 75252 PARIS Cedex 05; email: \texttt{amaury.lambert@upmc.fr};
  url: \texttt{http://www.lpma-paris.fr/pageperso/amaury.lambert/}
 }


\begin{abstract} Ultrametric trees are trees whose leaves lie at the same distance from the root. They are used to model the genealogy of a population of particles co-existing at the same point in time.  
We show how the boundary of an ultrametric tree, like any compact ultrametric space, can be represented in a simple way via the so-called comb metric. We display a variety of examples of random combs and explain how they can be used in applications. In particular, we review some old and recent results regarding the genetic structure of the population when throwing neutral mutations on the skeleton of the tree. 
\end{abstract}

\begin{resume} Les arbres ultramétriques sont les arbres dont les feuilles se trouvent toutes à la même distance de la racine. Ces arbres sont utilisés pour modéliser la généalogie d'une population de particules qui co-existent à un temps donné. 
Nous montrons que la frontière d'un arbre ultramétrique, comme tout espace ultramétrique compact, peut être représentée simplement via ce que nous appelons la distance de peigne. Nous examinons plusieurs exemples de peignes aléatoires et nous expliquons comme ils peuvent être utilisés dans certaines applications. En particulier, nous évoquons quelques résultats anciens ou plus récents concernant la structure génétique de la population lorsque l'on jette des mutations neutres sur le squelette de l'arbre.

\end{resume}

\thanks{The author thanks the  \emph{Center for Interdisciplinary Research in Biology} (Coll\`ege de France) for funding.
}
\maketitle


\bigskip


\tableofcontents

\section{Introduction}

In this paper, we review some mathematical properties of random tree models bearing in mind potential applications to evolutionary biology. Trees are used in \emph{population genetics}, to trace the genealogy of a set of homologous genes, also of individuals sampled from asexual populations (in contrast to genealogies of individuals from sexual populations, which are called \emph{pedigrees} and are not trees); in \emph{phylogenetics}, to represent the ancestral relationships between species; in \emph{epidemiology}, to model both the history of transmissions in epidemics and the genetic relationships between pathogenic strains. When these entities (genes, individuals, species, patients, pathogens) are sampled at the same point in time, simply called the present time, these trees are said \emph{ultrametric}, which means that all the \emph{leaves} of the tree lie at the same graph distance from the root. Mathematically, what actually is ultrametric w.r.t. the graph distance is the set of leaves of the tree, called its \emph{boundary}. \\

From a theoretical point of view, a tree starts from one particle (coinciding with the root of the tree) and is generated by series of replication events (birth, speciation, transmission, division), which produce the so-called \emph{branching points} of the tree (points whose complementary has at least three connected components) and of termination events (death, extinction, recovery, apoptosis), which produce the leaves, or \emph{tips}, of the tree (points whose complement is connected). Probabilistic models for these processes abound \cite{Lam08, Lam17}: branching processes, birth-death processes, Wright-Fisher and Moran model, lookdown process... The genealogy of particles present at time $t$ is the subtree spanned by all points at the same distance $t$ from the root. It is the ultrametric tree we have introduced in the previous paragraph, called the \emph{reduced tree} in probability and the \emph{coalescent tree} in population genetics. It is also the ball of radius $t$ centered at the root.\\

Seen from an empirical point of view, reduced trees are also called \emph{reconstructed trees}. Indeed, the tree itself is not available as data, and so has to be inferred (that is, reconstructed), typically from \emph{multiple sequence alignments}. This can be done because there exist measures of \emph{genetic distance} (i.e., dissimilarity between two aligned genetic sequences) which are good proxies of their \emph{genealogical distance} (i.e., graph distance, or twice the time since most recent common ancestor). In most organisms, most mutations (said neutral) occur along the lineages of the tree at a more or less constant pace, hence the name of \emph{molecular clock}. Modeling mutational processes by Poisson processes with possibly variable mutation rates in time and across genes then provides mathematical relations between genealogical distance and genetic distance, which can in turn be used to infer the tree from the sequences. Statistical methods for the inference of trees from sequences form a scientific field in its own right and will not be considered here any longer. Note that the empirical trees represented in Figures \ref{fig:genetree} and \ref{fig:phylodynamics} are not ultrametric even though the sequences used to infer them do co-exist at present time. This is because evolutionary trees are often represented with genetic distances rather than genealogical distances. \\

The study of phylogenetic trees has fueled much mathematical research in graph theory, geometry and probability \cite{DMT96, EvaBook, SSBook}. Here, we review some results, essentially those recently obtained by the author and his co-authors, around the study of ultrametric trees, initially motivated by two questions: the inference of the process most likely to have generated a given reduced tree (phylogenetics, epidemiology); the neutral genetic composition expected to be observed in a branching population (population genetics).\\

We start with the general definition of real tree as a metric space, then we introduce the reduced tree (root-centered ball) and its boundary (root-centered sphere). We then show how the boundary of an ultrametric tree, like any compact ultrametric space, can be represented in a simple way via the so-called comb metric. We display a variety of examples of deterministic and random combs, the infinite $p$-ary tree, the Kingman comb and comb-based exchangeable combs in general, the boundary of a branching process and coalescent point processes in general. In the last section, we review some old and recent results regarding the genetic structure of the population when throwing neutral mutations on the skeleton of the tree.

\afterpage{
\begin{figure}[!ht]
\includegraphics[width=7cm]{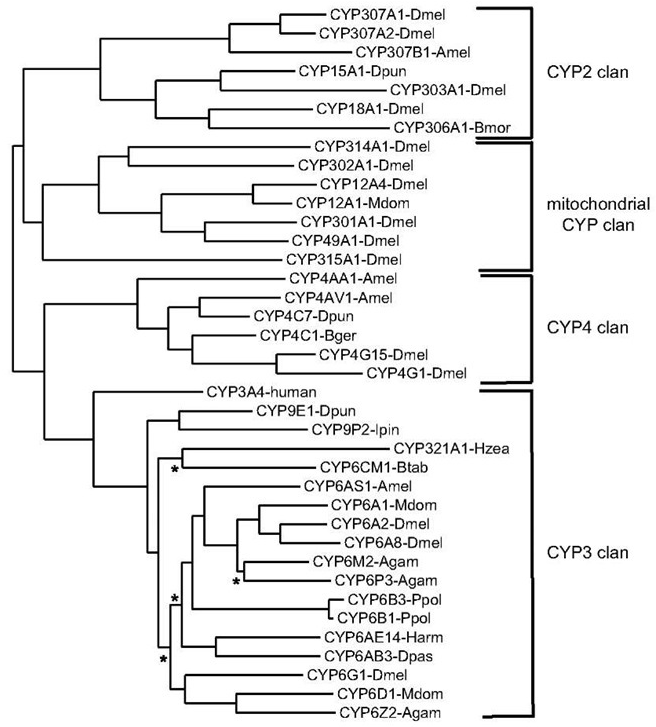}
\includegraphics[width=8cm]{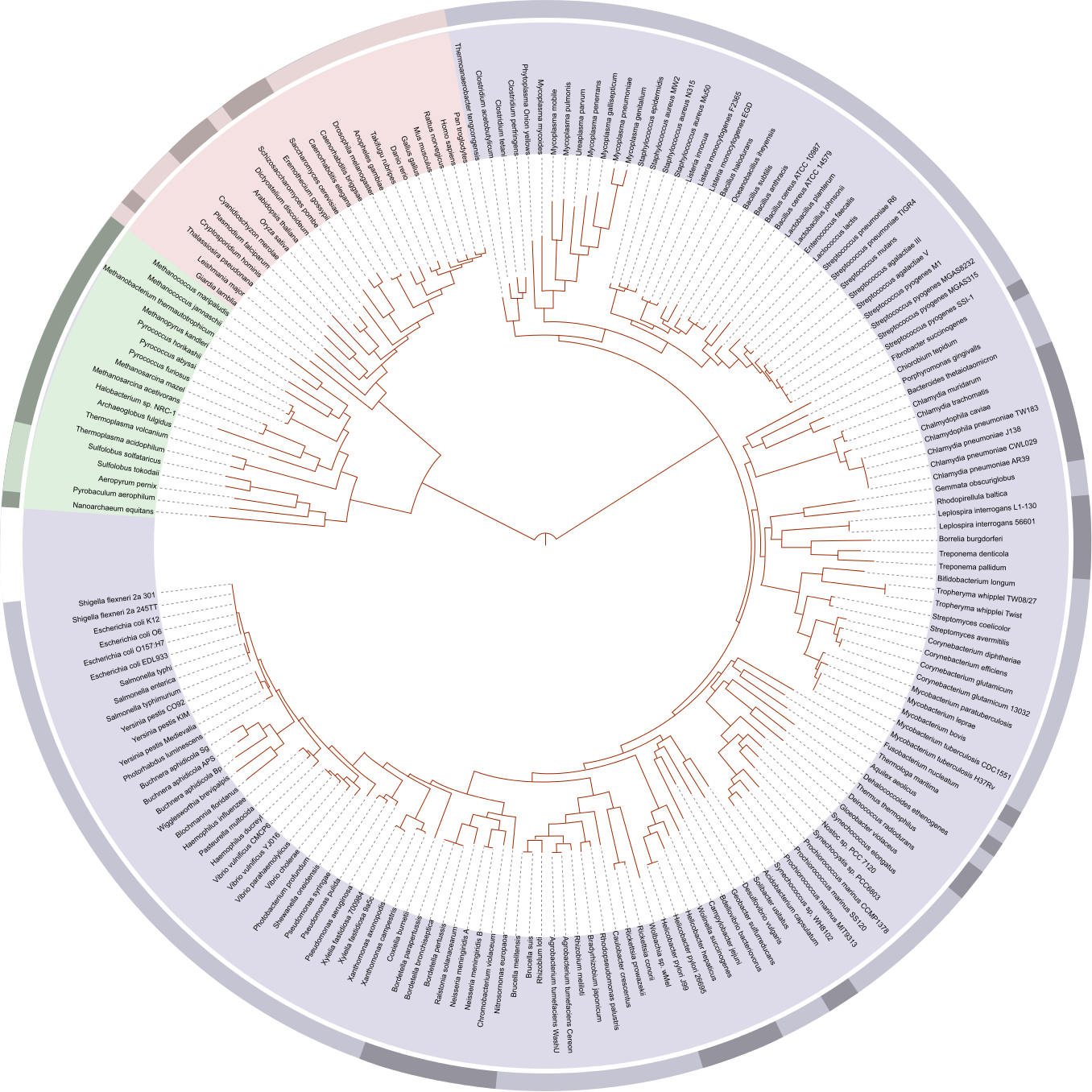}
\caption{Left panel: A gene tree showing the evolution of the CYP gene family in insects. 
Right panel: A phylogenetic tree of contemporary species, with different colors for each of the three domains: Bacteria (grey), Archaea (green) and Eukarya (pink)\protect\footnotemark.
}
\label{fig:genetree}
\end{figure}
\footnotetext{Sources: \url{http://what-when-how.com/insect-molecular-biology-and-biochemistry/insect-cyp-genes-and-p450-enzymes-part-2/} (left panel) and \url{https://en.wikipedia.org/wiki/Tree_of_life_(biology)} (right panel).}
}

\afterpage{
\begin{figure}[!ht]
\unitlength 1.8mm 
\linethickness{0.2pt}
\begin{picture}(54,28)(0,0)
\put(4,3){\line(1,0){50}}
\put(9,15){\line(0,1){6}}
\put(13,19){\line(0,1){8}}
\put(17,12){\line(0,1){5}}
\put(21,15){\line(0,1){5}}
\put(25,17){\line(0,1){11}}
\put(29,10){\line(0,1){8}}
\put(33,16){\line(0,1){9}}
\put(41,14){\line(0,1){5}}
\put(45,16){\line(0,1){4}}
\put(17,17){\circle*{.8}}
\put(9,21){\circle*{.8}}
\put(37,20){\circle*{.8}}
\put(8.93,14.93){\line(-1,0){.8}}
\put(7.33,14.93){\line(-1,0){.8}}
\put(5.73,14.93){\line(-1,0){.8}}
\put(12.93,18.93){\line(-1,0){.8}}
\put(11.33,18.93){\line(-1,0){.8}}
\put(9.73,18.93){\line(-1,0){.8}}
\put(16.93,11.93){\line(-1,0){.9231}}
\put(15.084,11.93){\line(-1,0){.9231}}
\put(13.237,11.93){\line(-1,0){.9231}}
\put(11.391,11.93){\line(-1,0){.9231}}
\put(9.545,11.93){\line(-1,0){.9231}}
\put(7.699,11.93){\line(-1,0){.9231}}
\put(5.853,11.93){\line(-1,0){.9231}}
\put(20.93,14.93){\line(-1,0){.8}}
\put(19.33,14.93){\line(-1,0){.8}}
\put(17.73,14.93){\line(-1,0){.8}}
\put(24.93,16.93){\line(-1,0){.8}}
\put(23.33,16.93){\line(-1,0){.8}}
\put(21.73,16.93){\line(-1,0){.8}}
\put(28.93,9.93){\line(-1,0){.96}}
\put(27.01,9.93){\line(-1,0){.96}}
\put(25.09,9.93){\line(-1,0){.96}}
\put(23.17,9.93){\line(-1,0){.96}}
\put(21.25,9.93){\line(-1,0){.96}}
\put(19.33,9.93){\line(-1,0){.96}}
\put(17.41,9.93){\line(-1,0){.96}}
\put(15.49,9.93){\line(-1,0){.96}}
\put(13.57,9.93){\line(-1,0){.96}}
\put(11.65,9.93){\line(-1,0){.96}}
\put(9.73,9.93){\line(-1,0){.96}}
\put(7.81,9.93){\line(-1,0){.96}}
\put(5.89,9.93){\line(-1,0){.96}}
\put(32.93,15.93){\line(-1,0){.8}}
\put(31.33,15.93){\line(-1,0){.8}}
\put(29.73,15.93){\line(-1,0){.8}}
\put(36.93,11.93){\line(-1,0){.8889}}
\put(35.152,11.93){\line(-1,0){.8889}}
\put(33.374,11.93){\line(-1,0){.8889}}
\put(31.596,11.93){\line(-1,0){.8889}}
\put(29.819,11.93){\line(-1,0){.8889}}
\put(40.93,13.93){\line(-1,0){.8}}
\put(39.33,13.93){\line(-1,0){.8}}
\put(37.73,13.93){\line(-1,0){.8}}
\put(44.93,15.93){\line(-1,0){.8}}
\put(43.33,15.93){\line(-1,0){.8}}
\put(41.73,15.93){\line(-1,0){.8}}
\multiput(3.93,22.93)(.979167,0){49}{{\rule{.4pt}{.4pt}}}
\put(3,23.375){\makebox(0,0)[cc]{\scriptsize $t$}}
\put(2.875,3.375){\makebox(0,0)[cc]{\scriptsize $0$}}
\put(37,12){\line(0,1){8}}
\put(5,18){\line(0,-1){15}}
\end{picture}
\includegraphics[width=7cm]{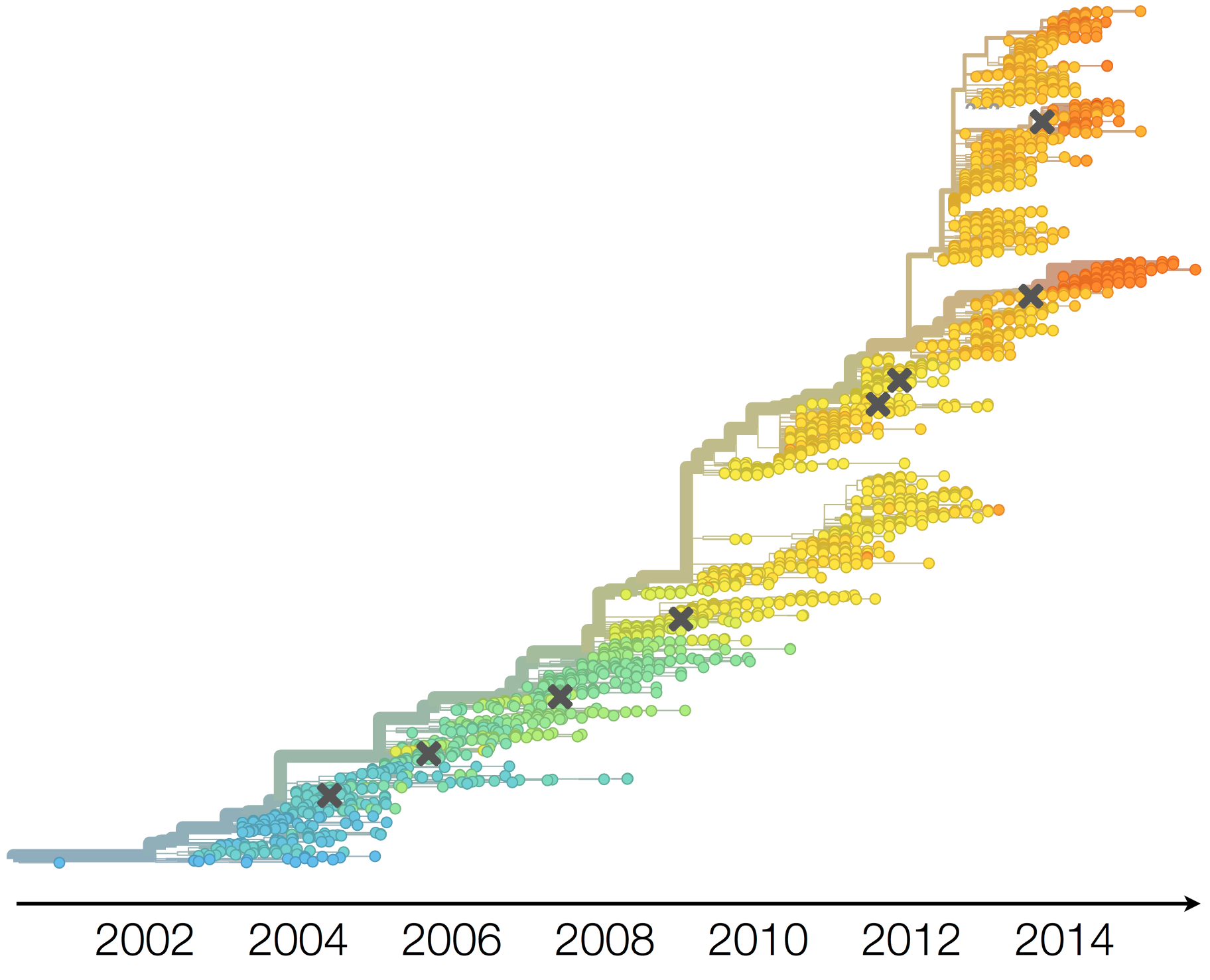}
\caption{Left panel: A hypothetical transmission tree, with dots showing recovery events (exit from the epidemic). Right panel: A phylogenetic tree of H3N2 strains of influenza\protect\footnotemark.}
\label{fig:phylodynamics}
\end{figure}
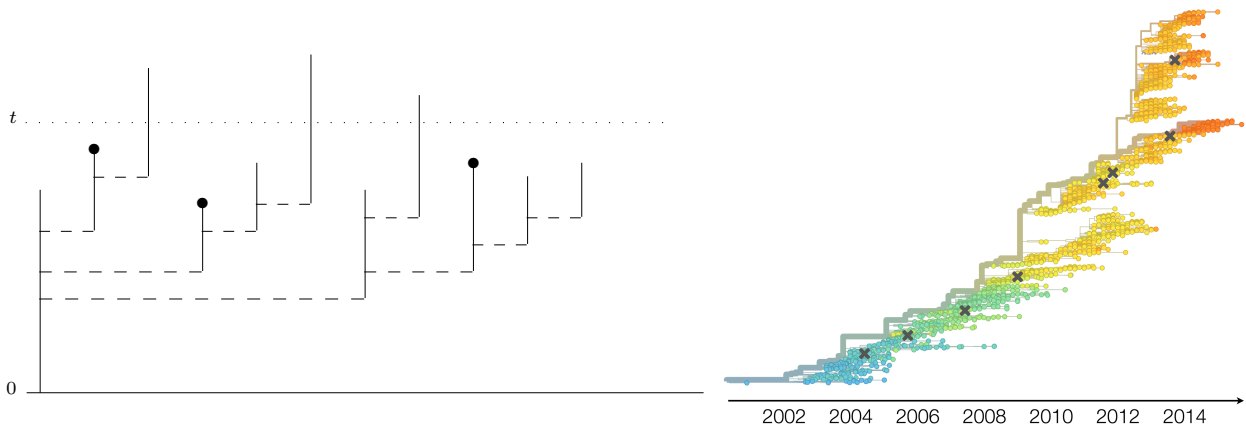
\footnotetext{Source (right panel): \url{http://bedford.io/talks/real-time-tracking-vidd/}}
}

\section{Real Trees, Ultrametric Trees, Combs}

\subsection{The real tree}

\begin{dfntn} A real tree, or $\R$-tree,  is a complete metric space $(\tr,d)$ satisfying
\begin{itemize}
\item[(A)] Uniqueness of geodesics. For any $x,y\in \tr$, there is a unique isometric map $\phi_{x,y}:[0,d(x,y)]\to\tr$ such that $\phi_{x,y}(0)= 0$ and $\phi_{x,y}(d(x,y))=y$. 

The geodesic $\phi_{x,y}([0, d(x,y)])$, also called \emph{arc}, is denoted $\llbracket x,y\rrbracket$. 
\item[(B)] No loop. For any continuous, injective map $\psi:[0,1]\to\tr$, $\psi([0,1])= \llbracket \psi(0),\psi(1)\rrbracket$.
\end{itemize}
The \emph{root} of an $\R$-tree $\tr$ is a distinguished element of $\tr$ denoted $\rho$.
\end{dfntn}

\begin{thrm}[Four points condition] The metric space $(\tr,d)$ is a real tree if it is complete, path-connected and satisfies for any $x_1, x_2, x_3, x_4\in\tr$
$$
d(x_1,x_2)+d(x_3,x_4)\le \max\{d(x_1,x_3)+d(x_2,x_4), d(x_1,x_4)+d(x_2,x_3)\}
$$
\end{thrm}
\noindent
For major references on this topic, see \cite{DMT96, EvaBook}.
\begin{dfntn}
For any $x\in\tr$, the \emph{multiplicity}, or \emph{degree} of $x$ denotes the number of connected components of $\tr\setminus\{x\}$. If 
$m(x)=1$, then $x$ is called a \emph{leaf} or a \emph{tip}, and if $m(x)\ge 3$, then $x$ is called a \emph{branching point}.
\end{dfntn}
We will further need the following notation and terminology.
\begin{itemize}
\item \textbf{Mrca}. For any $x,y\in\tr$ the \emph{most recent common ancestor} (in short mrca) of $x$ and $y$, denoted $x\wedge y$, is the unique $z\in \tr$ such that $\llbracket\rho,x\rrbracket \cap \llbracket\rho,y\rrbracket= \llbracket\rho,z\rrbracket$.
\item \textbf{Partial order}. For any $x,y\in\tr$, $y$ is said to \emph{descend} from $x$, and then $x$ is called an \emph{ancestor} of $y$ if $x\in\llbracket\rho,y\rrbracket$, and this is denoted $x\preceq y$. 
\item \textbf{Length measure}. Whenever $\tr$ is locally compact, there is a unique measure
 $\lambda$ on the Borel $\sigma$-field of $\tr$, called \emph{length measure}, such that for any $x,y\in\tr$, $\lambda(\llbracket x,y\rrbracket) = d(x,y)$ (see Section 4.3.5 in \cite{EvaBook}).
\item \textbf{Reduced tree}. For a real tree $\tr$ and a fixed real number $T>0$, the so-called \emph{reduced tree} at height $T$ is the tree spanned by points at distance $T$ from the root, i.e.
$$
\{y\in \tr: \exists x\in \tr, y\preceq x, d(\rho,x)=T\}.
$$
\end{itemize}
The topology of the reduced tree can be understood from the topology of the \emph{sphere} of $\tr$ with center $\rho$ and radius $T>0$
$$
\trT:= \{x\in\tr : d(\rho,x)=T\},
$$
Note that by the four-points condition, for any $x,y,z\in\trT$,
\begin{multline*}
T+d(x,z)=
d(\rho,y)+d(x,z)\le \max\{d(\rho,x)+d(y,z), d(\rho,z)+d(y,x)\} =\max \{T+d(y,z), T+d(y,x)\},
\end{multline*}
which yields $d(x,z)\le \max \{d(y,z), d(y,x)\}$, that is the metric induced by $d$ on $\trT$ is \emph{ultrametric}. 
From now on, we assume that $(\tr, d)$ is \textbf{locally compact}, so that \textbf{$(\trT,d)$ is a compact ultrametric space} (by application of the Hopf-Rinow theorem, since a real tree is a length-metric space). We will see in the next section that any compact ultrametric space can be represented by what we call a \emph{comb}.

\subsection{The comb metric}

Let $I$ be a compact interval and $f:I\to [0,\infty)$ such that for any $\varepsilon >0$, \textbf{$\{f\ge \varepsilon\}$ is finite}. 
For any $s,t\in I$, define $\dtf$ by
$$
\dtf(s,t) = 2
\max_{(s\wedge t, s\vee t)} f.
$$
It is clear that $\dtf$ is a pseudo-distance on $\{f=0\}$ and that it is ultrametric, i.e.
$$\dtf(r,t)\le \max\{\dtf(r,s),\dtf(s,t)\} \qquad r,s,t\in I.
$$
Let us assume additionally that  \textbf{$\{f\not=0\}$ is dense in $I$} for the usual topology, so that  $\dtf$ is a distance on $\{f=0\}$.
\begin{dfntn}
\label{dfntn:comb}
We call $f$ a \emph{comb-like function} or \emph{comb}, and $\dtf$ the \emph{comb metric} on $\{f=0\}$.
\end{dfntn}

\noindent
The space $(\{f=0\}, \dtf)$ is not complete in general. To make it complete, one has to distinguish for each point $t\in I$ between its \textbf{left face} $(t,l)$ and its \textbf{right face} $(t,r)$. The distance $\dtf$ is extended to the space $I\times\{l,r\}$ by the following definitions for $s< t\in I$
$$
\dtf((s,r),(t,l)) = 2\max_{(s, t)} f, \qquad \dtf((s,l),(t,l)) = 2\max_{[s, t)} f, 
$$
$$
\dtf((s,r),(t,r)) =2\max_{(s, t]} f,  \qquad \dtf((s,l),(t,r)) = 2\max_{[s, t]} f,
$$
and the symmetrized definitions for $s> t$. If $f(s)=0$, $\dtf((s,l), (s,r))=0$ so that $(s,l)$ and $(s,r)$ must be identified. It can be shown \cite{LUB17} that the associated \textbf{quotient space $(\bar{I}, \dtf)$ is a compact, ultrametric space} called \textbf{comb metric space}. 
Actually the converse also holds, as we will see with Theorem \ref{thrm:second sense}.\\

\begin{figure}[!ht]
\includegraphics[width=\textwidth]{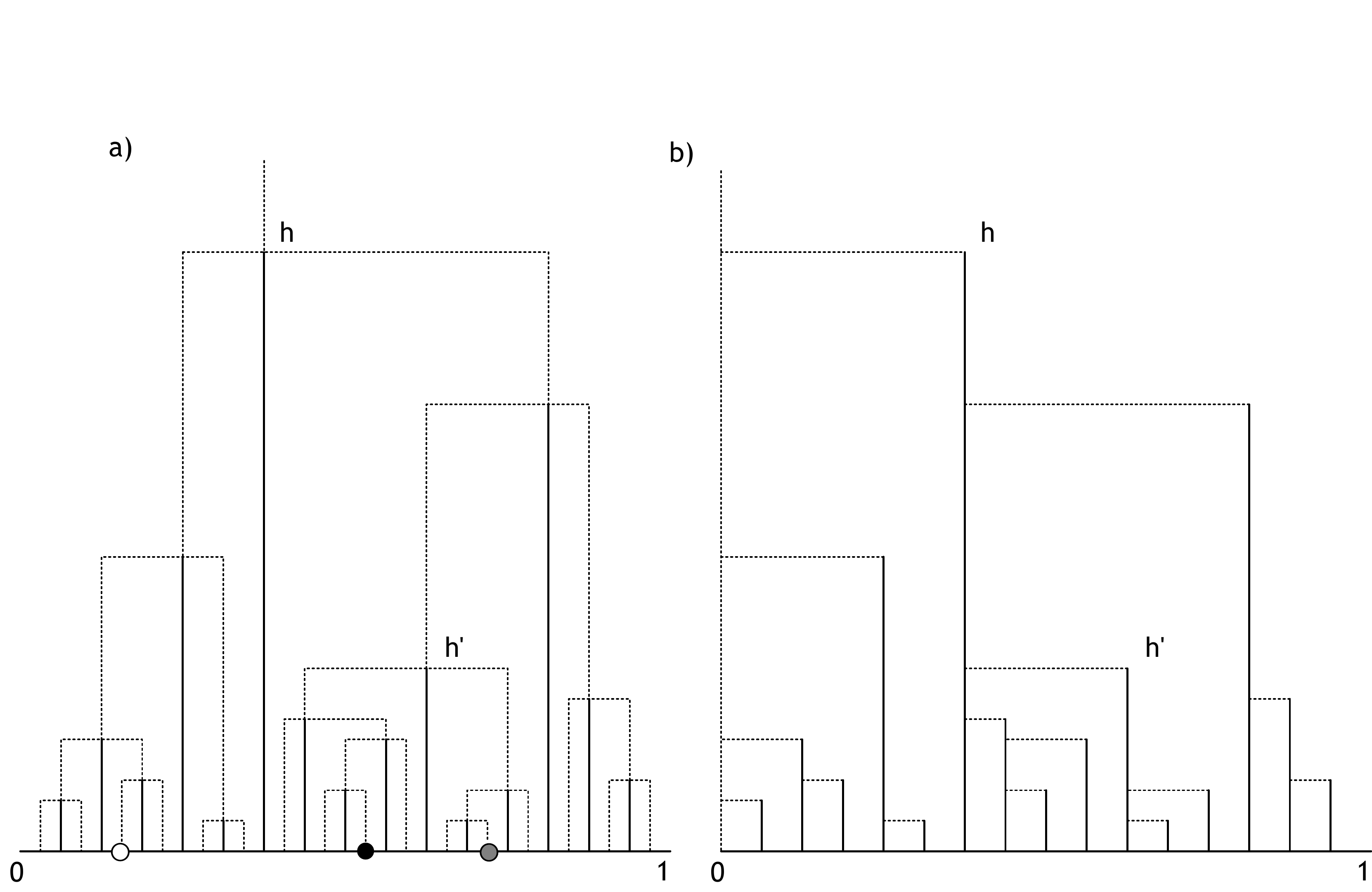}
\caption{A comb-like function on $[0,1]$ and two ways (a and b) of representing the associated ultrametric tree in dotted lines. a) The distance between the black dot and the grey dot is $h'$, whereas $h$ is the distance from either of these dots to the white dot; b) Graphical representation of the ultrametric tree constructed from the comb, introduced in Definition \ref{dfn:comb-tree}.}
\label{fig:comb}
\end{figure}

\noindent Before stating this theorem, we wish to construct the ultrametric tree hidden behind the comb metric space, as illustrated on Figure \ref{fig:comb}b.
\begin{dfntn} 
\label{dfn:comb-tree}
Let $f$ be a comb on the interval $I=[0,a]$ and $T>\max f$. We define 
$$
\text{Sk} := \{0\}\times(0,T] \cup \{(t,y) \in I\times(0,T]: \, f(t) > y\},
$$
endowed with the distance $\dtf$
\begin{gather*}
\dtf ((s,x), (t,y)) = \begin{cases}
			\;|\max_{(s,t]}f - x| +  |\max_{(s,t]}f - y| \quad & \text{if } s<t,\\
			\;|x-y| & \text{if } s=t.
		\end{cases}		
\end{gather*}
The tree $\tau_f(T)$ is defined as the \textbf{completion of $(\text{Sk}, \dtf)$}, so that in particular $\text{Sk}$ is its skeleton. We will always take the root of $\tau_f(T)$ equal to $\rho=(0,T)$.
For $t \in (0,T]$, we call the \textbf{lineage of $t$} the subset of the tree $L_t$ defined as the closure of the set
\[\{(s, x) \in \text{Sk}:\, s \leq t, \, x\ge\max_{(s,t]}f \}. \]
\end{dfntn}
\noindent
It can be shown \cite{DL17} that the boundary of $(\tau_f,\dtf)$ is indeed $(\bar{I}, \dtf)$, which explains why we keep the same notation for the two distances. Also that for each $t\in I$, $L_t=\llbracket \rho, \alpha_t\rrbracket$, where $\alpha_t\in \bar{I}$ is equal to $(t,r)$.

\begin{thrm}[\cite{LUB17}]
\label{thrm:second sense}
Any compact ultrametric space without isolated point is isometric to a comb metric space.
\end{thrm}
\noindent
Note that a comb metric space $(\bar{I}, \dtf)$ is naturally endowed with the \textbf{finite measure $\ell$} defined by 
$$
\label{dfn:measure}
\ell (A\times B) :=\text{Leb} (A),
$$
for any Borel set $A\subseteq I$ and $B\subseteq\{l,r\}$ (where $\text{Leb}$ denotes the Lebesgue measure, and it is important that $\{f\not=0\}$ is dense), which suggests that any compact ultrametric space can be equipped with a finite measure (isolated points can be treated separately). Actually, any compact ultrametric space can be endowed with a finite measure charging every ball with non-zero radius. One example of such a measure is the so-called \textbf{visibility measure} \cite{Lyo94}, as shown by the following argument, which is actually also used in the proof of the theorem. 

For any ultrametric space $(U,d)$ and any $r>0$, $U$ can be partitioned into balls of radius $r$, because the relation $\sim_r$ defined by $x\sim_r y \Leftrightarrow d(x,y)\le r$ is an equivalence relation. If in addition $U$ is assumed to be compact, the number $M_r$ of blocks in this partition has to be finite, and it is nondecreasing in $r$. The visibility measure is constructed by putting mass 1 on $U$, and recursively at each jump time of $M$ as $r$ decreases, by dividing the mass of each fragmenting block equally between its new sub-blocks. The comb can be constructed simultaneously with this recursive construction of the visibility measure, by mapping each ball with measure $m$ to an interval with length $m$, and putting `walls' between such intervals (the graph of the comb). See \cite{LUB17} for the mathematical details.

Of course, if $U$ is already endowed with a measure, the same construction can be done using the given measure. This is in particular the case when the ultrametric space is a sphere of a totally ordered, measured tree. It can then be shown \cite{LUB16} that there is a c\`adl\`ag function $h:\R\to\R$ with no negative jumps which codes for the tree in a sense that we specify hereafter.

\subsection{Sphere of a tree coded by a real function}

\begin{figure}[!ht]
\includegraphics[width=\textwidth]{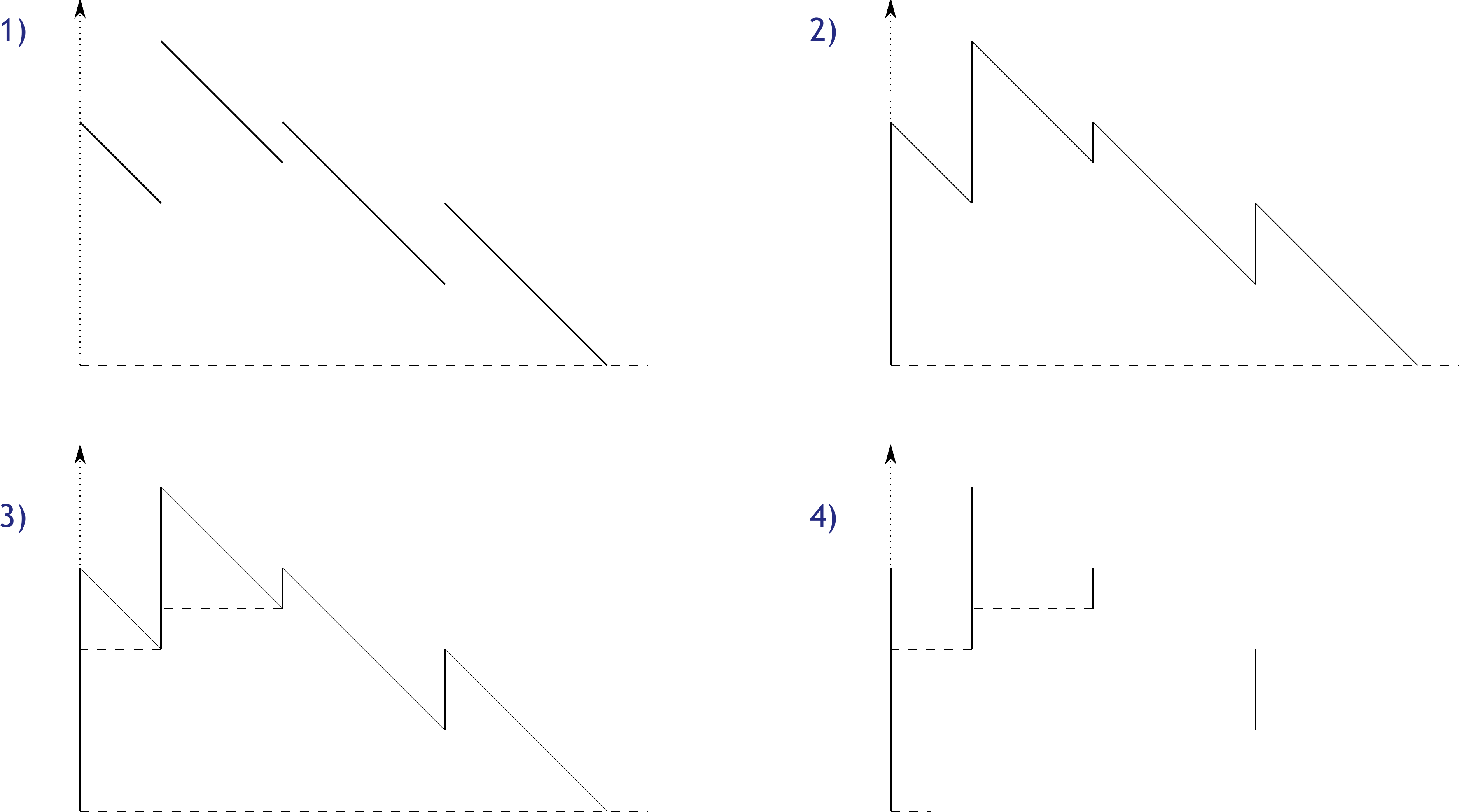}
\caption{How to uncover the tree coded in a càdlàg map which has a.e. negative derivative and positive jumps; 1) Start from the graph of this map; 2) Draw vertical solid lines in the place of jumps; 3) Report horizontal dashed lines from each edge bottom left to the rightmost solid point; 4) erase diagonal lines.}
\label{fig:jcp}
\end{figure}

Let $h:[0,\infty)\to[0,\infty)$ be càdlàg with no negative jumps and compact support. We are going to explain how $h$ codes for a real tree. 
Set $\sigma_h:=\sup\{t>0: h(t)\not=0\}$
and
$$
d_h(s,t):=h(s)+h(t) -2\inf_{[s\wedge t, s\vee t]}h.
$$
It is clear that $d_h$ is a pseudo-distance on $[0,\infty)$. Further let $\sim_h$ denote the equivalence relation on $[0,\infty)$
$$
s\sim_h t \Leftrightarrow d_h(s,t)=0 \Leftrightarrow h(s)=h(t)=\inf_{[s\wedge t, s\vee t]}h.
$$
\begin{thrm}
\label{thrm:coding}
Denote by $\tr_h$ the quotient space $[0,\sigma_h]|_{\sim_h}$. Then $(\tr_h,d_h)$ is a compact $\R$-tree. 
\end{thrm}
\noindent Figure \ref{fig:jcp} shows how to uncover the tree coded by a function, for a very simple example.
\noindent
From now on, let $p_h:[0,\sigma_h]\to \tr_h$ map any element of $[0,\sigma_h]$ to its equivalence class relative to $\sim_h$. Note that the tree $\tr_h$ is naturally endowed with a total order and a mass measure, as follows. 
\begin{itemize}
\item \textbf{Total order}. We define $\le_h$ as the \emph{order of first visits}, that is for any $x,y\in\tr_h$,
$$
x\le_h y \Leftrightarrow \inf p_h^{-1}(\{x\})\le \inf p_h^{-1}(\{y\}).
$$
\item \textbf{Mass measure}. The measure $\mu_h$ is defined as the \emph{push forward of Lebesgue measure} by $p_h$.
\end{itemize}
\noindent
Conversely, it can actually be shown \cite{Duq06, LUB16} that if $(\tr, d)$ is a compact $\R$-tree endowed with a total order $\le$ and a finite mass measure $\mu$ satisfying some consistency conditions, there is a unique càdlàg map $h$ called the \textbf{jumping contour process} of $\tr$ such that the tree $(\tr_h, d_h, \le_h, \mu_h)$ is isomorphic to $(\tr, d, \le, \mu)$ (follow the panels of Figure \ref{fig:jcp} in reverse order).\\

Now let us consider a compact real tree $\tr$ coded by a function $h$ (its jumping contour process if the function is not given \textit{a priori}) and let $T>0$ be such that the sphere $\trT$ is not empty. We know that $\trT$ is an ultrametric space, but we have no guarantee that $\mu_h$ charges $\trT$, so we will directly construct an isometric comb metric space, following \cite{LUB17}.

Assume that the sphere $\trT$ has no isolated point in itself. Then it can be shown that $\{h=T\}$ has no isolated point and has empty interior. So we can construct a \emph{local time} at level $T$ for $h$, that is a nondecreasing, continuous map $L:[0,\infty)\to [0,\infty)$ such that $L(0)=0$ and for any $s<t$
$$
L_t>L_s \Leftrightarrow (s,t)\cap\{h=T\}\not=\varnothing.
$$
Let $I=[0,L_{\sigma_h}]$. Let $J$ denote the right inverse of the local time 
$$
J_t := \inf\{s>0:L_s >t\}\qquad t\in I.
$$
Every jump time $s$ of $J$ corresponds to an interval $(g_s,d_s):=(J_{s-}, J_s)$ where $L$ is constant (to $s$) such that $h\not=T$ on $(g_s,d_s)$ and $h(g_s-)=h(g_s)= T = h(d_s)$ (because $h$ has no negative jump). In particular, $\trT$ is the closure of the range of $p_h\circ J$. Now for any $s<t$ such that $h(s)=h(t)=T$,
$$
d_h(s,t)  = h(s)+h(t) -2\inf_{[s, t]}h = 2(T-\inf_{[s,t]}h).
$$
so for any $s<t$ in $I$, the distance between the two points $p_h(J_s)$ and $p_h(J_t)$ of $\trT$ is $2\max_{[s,t]} f$, where
$$
f(s):=
\left\{
\begin{array}{cl}
T-\inf_{[g_s, d_s]}h &\text{ if }\Delta J_s \not=0\\
0 &\text{ otherwise.}
\end{array}
\right.
$$
This indicates that $\trT$ should be represented by the comb metric space associated to the comb $f$ just defined.
Note that $f$ is a comb with values equal to the \textbf{depths of the excursions of $h$ away from $T$}.
Further define
\Fonction{\bar\theta}{(\bar I,\dtf)}{(\trT, d)}{s'}{
\left\{
\begin{array}{cl}
p_h(J_{s-})& \text{ if $s'=(s,l)$ and $f(s)\not=0$}\\
 p_h(J_s)& \text{ otherwise.} 
\end{array}
\right.
}
\begin{thrm}[\cite{LUB17}] 
\label{thrm:comb-reducedtree}
The map $\bar\theta:(\bar I, \dtf)\to(\trT,d_h)$ is a global isometry preserving the order and mapping the Lebesgue measure to the push forward of the measure $dL$ by $p_h$.
\end{thrm}

\subsection{The boundary of the infinite $p$-ary tree}
\begin{figure}[!ht]
\includegraphics[width=16cm]{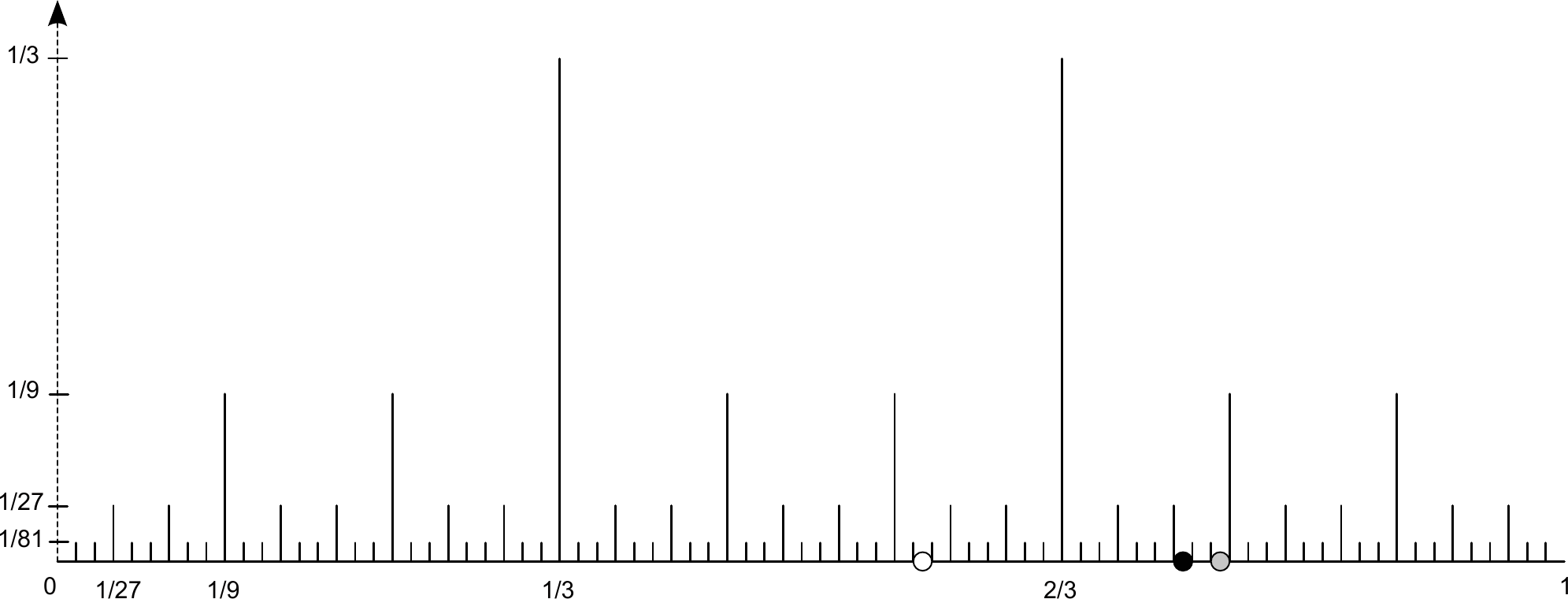}
\caption{Triadic comb. The white dot represents a sequence starting with $(1,2,0,1,\ldots)$, the black dot represents a sequence starting with $(2,0,2,0,\ldots)$ and the grey dot a sequence starting with $(2,0,2,2,\ldots)$. The distances white/black and white/grey are both equal to $(1/3)^1$, and the distance black/grey equals $(1/3)^4$.}
\label{fig:TriadicComb}
\end{figure}

The fundamental example of compact ultrametric space is the boundary of the infinite $p$-ary tree. This is the set $U_p$ of sequences $x=(x_n)_{n\ge 1}$ with values in $\{0,\ldots,p-1\}$ endowed with the distance $d_u$ defined for any pair $x=(x_n)$ and $y=(y_n)$ of elements of $U_p$ by
$$
d_u(x,y) = p^{-v_u(x-y)},
$$
where 
$$
v_u(x):=\min\{n\ge 1:x_n\not=0\}, 
$$
with the convention $p^{-v_u(0)}=0$.  
It is known that $(U_p,d_u)$ is a compact ultrametric space. The distance $d_u$ is actually the graph distance associated to the tree where each edge between generation $n-1$ and generation $n$ has length equal to $p^{-n}$. We show how to construct explicitly the isometry of the theorem between $(U_p,d_u)$ and a comb metric space which can intuitively be guessed from the previous remark, as shown on Figure \ref{fig:TriadicComb}.

\begin{figure}[!ht]
\includegraphics[width=16cm]{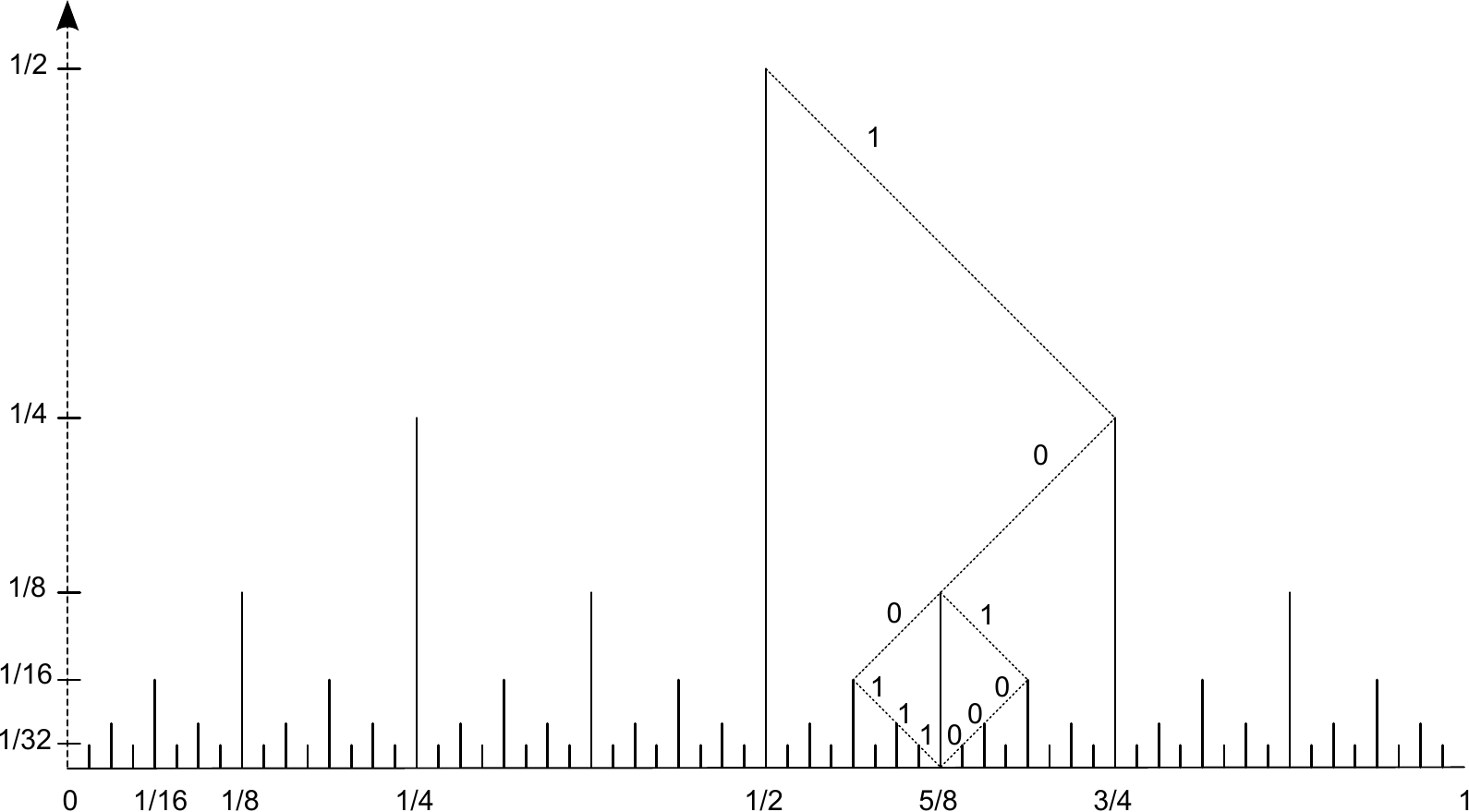}
\caption{Dyadic comb. The dashed lines represent two infinite sequences of 0's and 1's: on the right, the sequence $x=(1,0,1,0,0,0,\ldots)$ and on the left the sequence $\hat x=(1,0,0,1,1,1,\ldots)$. They correspond to the right and left faces of 5/8, i.e., the dyadic decompositions of 5/8 stationary at 0 and at 1, respectively. They lie at distance $(1/2)^3$ from each other.}
\label{fig:DyadicComb}
\end{figure}

\noindent
Recall that $t\in [0,1]$ is called a \emph{$p$-adic number} if it has two distinct $p$-adic decompositions. We denote by $\phi_l(t)$ its $p$-adic decomposition stationary at $p-1$ and by $\phi_r(t)$ its $p$-adic decomposition stationary at $0$. When $t$ is not $p$-adic, its $p$-adic decomposition is simply denoted $\phi(t)$.
Now for any $x\in U_p$, set $w(x):=\max\{n\ge 1 :x_{n}\not=0\}$,
and for any $t\in[0,1]$, define
$$
F_p(t):=
\left\{
\begin{array}{cl}
p^{-w(\phi_r(t))} &\mbox{ if $t$ is $p$-adic}  \\
0 &\mbox{ otherwise.}
\end{array}
\right.
$$
The function $F_p$ is a comb called the \textbf{$p$-adic comb} and it is not difficult to see that $\dtfp$ defines a comb metric space $(\bar{I},\dtfp)$. Then we define $\bar{\phi}:(\bar{I},\dtfp)\to(U_p,d_u)$ that maps every $t'\in[0,1]\times\{l,r\}$ to the $p$-adic decomposition of $t$, except when $t$ is $p$-adic which then has two faces, each mapped to a distinct $p$-adic decomposition as follows.
$$
\bar{\phi}(t'):=
\left\{
\begin{array}{cl}
\phi(t) &\mbox{ if } t'=t \in \{F_p=0\}\\
\phi_l(t) &\mbox{ if } t' =(t,l)\\
\phi_r(t)& \mbox{ if } t' =(t,r).
\end{array}
\right.
$$
It can be seen that $\bar{\phi}$ is a global isometry between $(\bar{I},\dtfp)$ and $(U_p,d_u)$ conserving the measures, where $U_p$ is endowed with the visibility measure $\mu_p$ that puts mass $p^{-n}$ on any ball with radius $p^{-n}$. The measure $\mu_p$ is also the law of a sequence of i.i.d. random variables uniformly distributed in $\{0,\ldots, p-1\}$. See Figure \ref{fig:DyadicComb} for an illustration of the \textbf{dyadic comb} showing the left and right faces of the same dyadic number.

\section{Random Ultrametric Trees}

\subsection{The comb-based exchangeable coalescent}
\label{subsec:CEC}

A first example of random ultrametric tree is the following.
Start with a comb $f$ on $[0,1]$ and an independent sequence $(V_i)$ of independent and identically distributed (i.i.d.) random variables uniform in $(0,1)$. Note that a.s. $f(V_i)=0$ for all $i$, so we can define the ultrametric distance $\delta_f$ on $\N$ by
$$
\delta_f(i,j):= \dtf(V_i, V_j).
$$ 
Now for any $t>0$, define  the partition $R_f(t)$ on $\N$ induced by the equivalence relation $\sim_t$
$$
i\sim_t j \Leftrightarrow \delta_f(i, j) \le t.
$$
Note that  $R_f(0+)$ is the finest partition (all singletons) because $\{f\not=0\}$ is assumed dense in $[0,1]$.
\begin{dfntn}
The \emph{comb-based coalescent process} $(R_f(t);t>0)$ is an \emph{exchangeable coalescent process}, in the sense that its law is invariant under permutations of $\N$.
\end{dfntn}
\noindent
The converse statement is given in the next proposition.
\begin{prpstn}
Let $(R(t);t>0)$ be a random exchangeable coalescent process such that for each $t>0$, $R(t)$ has a finite number of blocks and no singleton. Then there is a random comb $f$ such that the comb-based coalescent process $R_f$ is equal in distribution to $R$.
\end{prpstn}
\begin{proof}
We can define the distance $\delta$ on $\N$ by
$$
\delta(i,j):= \inf\{t\ge 0: \mbox{ $i$ and $j$ are in the same block of }R(t)\}.
$$ 
Then it is straightforward that $\delta$ is ultrametric. Indeed, for any integers $i, j, k$, and for any $t\ge 0$, if $\delta(i,j)\le t$ and $\delta(j,k)\le t$, then $i$ and $j$ on the one hand, and $j$ and $k$ on the other hand, are in the same block of $R(t)$, so that $i$ and $k$ are in the same block of $R(t)$. This shows that $\delta(i,k)\le t$, so that $\delta(i,j)\le \max\{\delta(i,j),\delta(j,k)\}$. Now by de Finetti's theorem, each ball $A$ of the ultrametric space $(\N,\delta)$ has an asymptotic frequency $m(A)$ defined by
$$
m(A):=\lim_{n\uparrow\infty}\frac{\#A\cap\{1,\ldots, n\}}n.
$$ 
Unfortunately, $m$ is not a measure (it does not satisfy Caratheodory's property) and $(\N,\delta)$ is not compact. But because we assume that for all $t>0$, $R(t)$ has a finite number of blocks and no singleton, we can use the same procedure as outlined after the statement of Theorem \ref{thrm:second sense} and construct a comb $f$ on $I=[0,1]$ and a map $\phi$ between balls of $(\N,\delta)$ and balls of $(I,f)$, which conserves radii, measures and partial order (inclusion). Now let $(V_i)$ be i.i.d. uniform random variables independent of the comb $(I,f)$ and of the map $\phi$ and define $\delta_f$ by  $\delta_f(i,j)= \dtf(V_i,V_j)$.
By construction, the processes of ranked frequencies of blocks of $R$ and of $R_f$ are a.s. equal. Kingman's representation theorem of exchangeable partitions ensures that the law of $R(\varepsilon)$ can be obtained from its frequencies by a paintbox process, so that $R_f(\varepsilon)$ has the same law as $R(\varepsilon)$. By monotonicity, this implies that $(R_f(t);t\ge\varepsilon)$ and $(R(t);t\ge\varepsilon)$ are equally distributed. Since $\varepsilon$ is arbitrary, the result is proved.
\end{proof}
\noindent
The archetypal example of comb-based coalescent is the following.
\begin{dfntn}
\label{dfn:kingman-comb}
Define the random comb $f$ on $I=[0,1]$ by
$$
f = \sum_{j\ge 1} \tau_j \indicbis{U_j},
$$
where the families of r.v. $(U_j)$ and $(\tau_j)$ are independent, the $(U_j)$ are i.i.d. uniform on $(0,1)$ and $\tau_j=\sum_{k\ge j+1} e_k$, where $e_k$ are independent exponential r.v. with parameter $k(k-1)/2$. Then the comb-based coalescent $R_f$ has the same distribution as the \textbf{Kingman coalescent} \cite{Kin82}. We call $f$ the \textbf{Kingman comb}.
\end{dfntn}
\noindent
In the rest of this section, we assume that $\tr$ denotes a \textbf{binary} $\R$-tree, that is, $m(x)\le 3$ for all $x\in\tr$. See  \cite{LP13} for extensions to random real trees with arbitrarily large degree.  

\subsection{The coalescent point process}

Recall from Theorem \ref{thrm:coding} that any càdlàg map $h$ codes for a compact $\R$-tree denoted $\tr_h$. Let us call \textbf{Brownian tree} the tree coded by a Brownian excursion conditioned to have height larger than $T$. The Brownian excursion has a local time at level $T$, which allows one to construct as in Theorem \ref{thrm:comb-reducedtree} the comb giving the metric of the reduced tree at level $T$.  This comb is a `list', in the plane order, of the depths of excursions of the contour away from $T$. These depths form a Poisson point process, which motivates the following definition.

\begin{dfntn}
Let $\nu$ be a $\sigma$-finite measure on $(0,\infty)$ such that $\nu([\varepsilon,\infty))<\infty$ for all $\varepsilon>0$. Let $\mc$ be a Poisson point process on $(0,\infty)^2$ with intensity $\mbox{Leb}\,\otimes\,\nu$ and denote by $(S_i, H_i)_i$ its atoms. Finally, let $(D,H)$ denote the first (i.e., with smallest first component) atom  such that $H>T$.
We will say that the random comb metric space associated with the comb $\sum_{i:S_i<D}H_i\indicbis{S_i}$ is a \textbf{coalescent point process} with height $T$ and intensity measure $\nu$.
\end{dfntn}
\begin{prpstn}[\cite{Pop04, AP05}]
\label{prop:brownian-cpp}
The sphere $\trT$ of the Brownian tree $\tr$ is isometric to a coalescent point process with height $T$ and intensity measure $\nu_0$, where $\nu_0$ is the push forward of the Brownian excursion measure by the function which maps an excursion to its depth, i.e.
$$
\nu_0(dx) = \frac{dx}{2 x^2} 
$$
We will call \textbf{Brownian CPP} a CPP with intensity measure $\nu_0$ (or a multiple of it).
\end{prpstn}
\noindent
More generally, a (root-centered) \textbf{sphere of any real tree whose contour process is strongly Markovian is isometric to a coalescent point process (CPP)}. This is the case of \textbf{splitting trees}, which are the trees generated by binary branching processes where particles give birth at constant rate during a lifetime which follows an arbitrary distribution \cite{GK97}. A splitting tree is actually isometric to a tree coded by a Lévy process with finite variation \cite{Lam10}.

More generally, we say that a random binary $\R$-tree satisfies the \textbf{splitting property} if for any $x\in \tr$, the subtrees rooted at the branching points of the segment $\llbracket\rho, x\rrbracket$ form a Poisson point process on $\llbracket\rho, x\rrbracket\times \mathcal E$, where $\mathcal E$ denotes the space of locally compact $\R$-trees. For example, the tree coded by a Brownian excursion satisfies the splitting property. We have proved that a tree satisfying the splitting property, again called splitting tree, is isometric to a tree coded by a Lévy process with possibly infinite variation \cite{LUB16}, so that its sphere is again isometric to a CPP.

From a practical point of view, an interesting question is to characterize the intensity measure $\nu$ of the CPP. In the case of strongly Markovian contour processes, \textbf{$\nu$ is simply the push forward of the excursion measure away from $T$ of the contour process, by the function which maps an excursion to its depth}, and a lot is known on this measure in the Lévy case (see \cite{Lam10} for details). There is also a whole class of random trees whose spheres are CPPs with finite intensity measure (derogating to our definition of combs having dense support) and which have applications in evolutionary biology.

Consider a population where all individuals live and reproduce independently, and each individual is endowed with a trait (some random character living in $\R$ for simplicity) that evolves through time according to independent copies of the same, possibly time-inhomogeneous, Markov process. Further assume what follows.
\begin{itemize}
\item This trait is non-heritable, in the sense that any individual born at time $t$ draws the value of her trait at birth from the same distribution, independently of her ancestors' histories;
\item All individuals give birth during their lifetime, according to a Poisson point process with intensity $\beta$, where $\beta$ is a diffuse Radon measure on $[0,\infty)$;
\item An individual holding trait $x$ at time $t$ dies at rate $d(t,x)$. 
\end{itemize}
\begin{thrm}[\cite{LS13}]
\label{thrm:generalmodel}
For a tree generated by the previously defined population model, starting with one individual at time 0 and conditional on having at least one alive individual at time $T$, the sphere of radius $T$ is isometric to a CPP with intensity measure $\nu$ given by
$$
W(t):=\frac{1}{\nu([t,\infty))} =\exp\left( \int_{T-t}^ T  (1-q(s))\,\beta(ds)\right)\qquad t\in[0,T],
$$
where $q(t)$ denotes the probability that an individual born at time $t$ has no descendants alive by time $T$. Note that the knowledge of $W(t)$ for $t>T$ is not needed.
In addition, if $g(t,\cdot)$ denotes the density of the death time of an individual born at time $t$, then $W$ is solution to the following integro-differential equation
$$
W'(t) = b(T-t)\,\left( W(t) - \int_0^t \ W(s)\,g(T-t,T-s) ds \right)\qquad t\ge 0,
$$
with initial condition $W(0)=1$.
\end{thrm}
The last statement can be used to compute the likelihood of a given reconstructed tree (cf. Introduction) to infer the model that most likely has generated this tree. Two beautiful biology papers applying this procedure are \cite{Sta11, MPP11}. More recent examples of applications can be found in Section 3.3 of the lecture notes \cite{Lam17}, see in particular \cite{EML14, LME15, LAS14, ALS15}.

\subsection{The boundary of a pure-birth process}
\label{subsec:yule}

Here, we consider a pure-birth process on $[0,\infty)$ with birth rate $\beta$, assumed to be a diffuse Radon measure on  $[0,\infty)$ such that $\beta([0,\infty))=\infty$. We will write $\underline\beta(t):=\beta([0,t])$.

The genealogy of such a process can be encoded by the infinite binary tree  $\mathcal{T}$ (finite sequences of 0's and 1's) endowed with the birth dates $\alpha (v)$ of each finite sequence $v$. If $u$ is a prefix of $v$, we write $u\preceq v$ and we say that $v$ is a \emph{descendant} of $u$.  
We will denote  $\omega(v)$ the date of death of $v$ and assume that $\omega(v) = \alpha(v0) =\alpha(v1)$. 
The tree $\mathcal{T}$ may then be equipped with a \textbf{measure $\mathscr{L}$ on its boundary $\partial\mathcal{T}= \{0,1\}^\N$} defined by
  \[
   \mathscr{L}\left (B_u\right ) := \lim _{t\uparrow \infty} \frac{N_u(t)}{e^{\underline\beta(t)}} \qquad  u\in \mathcal{T},
  \]
  where $B_u$ denotes the set of infinite sequences $v$ with prefix $u$ and $N_u(t)$ is the number of descendants of $u$ at time $t$ 
  \[N_u(t) := \#\{v\in \mathcal{T}, \; u\preceq v, \, \alpha(v) < t \leq \omega(v)\}.\]
The mass of $\mathscr{L}$ is an exponential r.v. with parameter 1.

\begin{thrm} \label{thm_yule}
Let $\varphi:[0,\infty)\to (0,T]$ be a decreasing bijection and let $\varphi(\tr)$ denote the tree obtained from the previous pure-birth tree $\tr$ by mapping all distances to the root by $\varphi$ and letting the measure $\mathscr{L}$ unchanged. Then $\varphi(\tr)$ is isometric to the tree $\tau_f$ introduced in Definition \ref{dfn:comb-tree} constructed from a CPP $f$ with intensity $\nu$ and height $T$, naturally measured by $\ell$ (defined page \pageref{dfn:measure}), where  
$$
\nu([t,\infty)) = \exp\left(\underline\beta\circ\varphi^{-1}(t) \right)\qquad t\in[0,T].
$$
Note that $\nu([T,\infty))=1$ due to the normalization of $\mathscr{L}$.

In particular, if we take $T=1$ and $\varphi$ given by  $\varphi(t) =e^{-\underline\beta(t)}$, then $\varphi(\tr)$ has the distribution of $\tau_f$, where $f$ is a CPP with height 1 and intensity $\frac{d x}{x^2}$.
\end{thrm}
\begin{proof}
We know from \cite{DL17} that $\varphi(\tr)$ is the genealogy of a reversed (i.e. time flows from $T$ to 0) pure-birth tree with birth rate $\beta'=\beta\circ\varphi^{-1}$. From the same paper, we know that a reversed pure-birth tree with birth rate $\beta'$ and measure $\mathscr{L}$ (invariant by the time-change) on its boundary is the tree $\tau_f$ associated with a CPP $f$ with intensity measure $\nu$ satisfying that the Laplace-Stieltjes measure associated with the increasing function $\ln(\overline \nu)$ is equal to $\beta'$, where we wrote $\overline\nu(t)= \nu([t,\infty))$. The proof follows from a simple calculation.    
\end{proof}
\noindent
Notice that we could have applied the same reasoning to a supercritical birth-death process conditioned to non-extinction, by considering the subtree spanned by its boundary, i.e. the tree of indefinitely surviving lineages, which has the law of the tree generated by the pure-birth process with birth rate $\tilde\beta(dt)=\beta(dt)\,(1-q(t))$, where $q(t)$ is the probability of survival of a single particle alive at time $t$. Then everything works as in Theorem \ref{thm_yule} provided that the counter $N_u(t)$ is restricted to the particles alive at $t$ with indefinitely surviving descendance, and that $\beta$ is replaced with $\tilde\beta$ in the displayed formula. 

Furthermore, we can recover the first part of Theorem \ref{thrm:generalmodel} by replacing in the last paragraph `infinite survival' with `survival up to $T$' and by taking $\varphi(t)=T-t$ in Theorem \ref{thm_yule}.

\subsection{Link between Kingman coalescent and CPP}
Recall the \textbf{Kingman comb} and the \textbf{Brownian CPP} defined respectively in Definition \ref{dfn:kingman-comb} and Proposition \ref{prop:brownian-cpp}. Both objects code for the genealogy of a large exchangeable population, but the Kingman coalescent is based on the assumption of a stationary population  with  constant  size  (total  size  constraint),  whereas  in  the  CPP  the  size of the population is fluctuating like a branching process, and its foundation time is fixed (time constraint). Following \cite{LS16}, we  show that one of the two is embedded in the other. In an exchangeable population with large constant size, the descendance of a small subpopulation is blind to the total size constraint and it is constant in expectation (see for example Theorem 1 in \cite{BL06}).  Our goal is to state a backward-in-time version of the last informal observation, namely `in a large stationary population with constant size, the genealogy of a subpopulation with recent MRCA is given by a CPP', and to derive some consequences of this fact.\\

Let $K:[0,1]\to[0,\infty)$ denote the Kingman comb and let $C:[0,\infty)\to[0,\infty)$ denote the CPP with intensity measure $\nu(dx)=2x^{-2}dx$.
Now for each $\varepsilon>0$, we define $S_\varepsilon$ the linear operator mapping each real function $f$ to the function $S_\varepsilon(f):t\mapsto \varepsilon^{-1}f\left(\varepsilon^{-1}t\right)$. The next statement ensures that if we zoom out on the Kingman comb restricted to a small interval, we end up with the Brownian CPP. 
\begin{prpstn}[\cite{LS16}]
As $\varepsilon \to 0$, the convergence $S_\varepsilon(K)\Longrightarrow C$ holds weakly in law on every compact set.
\end{prpstn}
\noindent
The genealogy generated by a comb restricted to a small interval focuses with high probability on parts of the ultrametric tree which are closely related (small genealogical distances). Now conditional on $K$, let us sample $n$ points conditioned to have a time to mrca smaller than $\varepsilon$, and consider the genealogy of the whole subtree spanned by this sample (quenched conditional sampling). In other words, we let $0=x^\varepsilon_0<x_1^\varepsilon<\cdots< x^\varepsilon_{N_\varepsilon}<1=x^\varepsilon_{N_\varepsilon+1}$ denote the ranked enumeration of the finite set $\{K>\varepsilon\}$ and we let $B_i^{\varepsilon}$ denote the comb restricted to the $i$-th interval of the subdivision ($1\le i\le N_\varepsilon+1$)
$$
B_i^\varepsilon (x)=K\left(x_{i-1}^\varepsilon+x\right)\qquad 0< x< l_i^\varepsilon,
$$
where $l_i^\varepsilon:= x_i^\varepsilon-x_{i-1}^\varepsilon$ is the length of the $i$-th interval.
Then for any measurable $f:\Omega\times(0,\infty)\to[0,\infty)$ (where $\Omega$ denotes the set of combs endowed with the vague topology), we define conditional on $K$ 
$$ 
\mu^\varepsilon_K(f):= \frac{\sum_{i=1}^{N_\varepsilon} (l^\varepsilon_i)^nf\left(B_i^\varepsilon,l_i^\varepsilon\right)}{\sum_{i=1}^{N_\varepsilon} (l^\varepsilon_i)^n}
$$
We naturally extend the definition of $S_\varepsilon$ to bivariate $f$, by $S'_\varepsilon(f):(\omega, l)\mapsto (S_\varepsilon(\omega),\varepsilon^{-1}l)$.
\begin{thrm}[\cite{LS16}]
For any continuous $f:\Omega\times(0,\infty)\to[0,\infty)$,
$$
\lim_{\varepsilon\downarrow 0} \E(\mu^\varepsilon_K(S_\varepsilon'(f))) = \E(f(M,L)),
$$
where $L$ is a Gamma ($n+1,2$) r.v. and $M$ is an independent CPP with intensity measure $\nu(dx)\,\indic{x<1}$ restricted to the interval $[0,L]$.
\end{thrm}
\noindent
In words, the previous statement ensures that after proper rescaling and averaging over the Kingman coalescent, the subtree spanned by a quenched conditional sample can be described in terms of a $n$-size-biased (i.e., biased by the $n$-th power of its size) Brownian CPP with height 1.

\section{Ultrametric Trees with Neutral Mutations}

\subsection{Throwing point mutations on a tree}

As explained in the Introduction, phylogenetic trees are inferred from genetic distances, due to the existence of a so-called molecular clock which regulates the pace at which new mutations appear on the lineages of the tree. This provides statistical relationships between genealogical distances and genetic distances, that allow biologists to infer the former from the latter. \\

So we consider a point measure $M$ on the skeleton of a real tree $(\tr,d)$, that we will call mutation point measure, whose atoms are viewed as mutation events. Assume that each point $x$ of the tree is further given a type, or \textbf{allele}, inherited from the most recent atom of $M$ on $\llbracket \rho, x\rrbracket$, that is, the point
$$
\sigma(x):=\arg\max\{d(\rho,y): y\text{ atom of } M, y\preceq x \},
$$
which is set equal to $\rho$ if the above set is empty. This assumption is known as the \textbf{infinitely-many allele model}. A point which carries the same allele as the root will be said \textbf{clonal}. The partition of the boundary into distinct alleles is the so-called \textbf{allelic partition}.\\

There are usually two ways of studying the allelic partition in the framework of combs. 

One possibility is to consider the allelic partition of a sample of size $n$. This can readily be done as in Subsection \ref{subsec:CEC}, by considering $n$ r.v. $V_1,\ldots,V_n$ i.i.d. uniform in the interval $I$ of definition of the comb, and associate each $i$ with the most recent atom of the lineage $\llbracket\rho, V_i\rrbracket$. This induces a partition of $\{1,\ldots, n\}$ which can be described by the so-called \textbf{allele frequency spectrum} $(A(1),\ldots, A(n))$, where $A(k)$ is the number of blocks of the partition with cardinality $k$.

A second possibility is to consider the allelic partition of the whole population. If the comb has finite support (finite number of teeth), then one can proceed as previously. Otherwise, the allele frequency spectrum has to be expressed thanks to the measure $\ell$ defined on the boundary, by defining the point measure on $(0,\infty)$
$$
A(dx):=\sum_{m} \delta_{\ell(R_m)}(dx),
$$ 
where the sum is taken over all atoms $m$ of $M$ and $R_m$ denotes the set of points $x$ in the boundary such that $\sigma(x)=m$.

In the context of the molecular clock, the most natural way of modeling mutation events is to use Poisson point processes. 
We have seen that a locally compact $\R$-tree $(\tr, d)$ has a length measure $\lambda$, so for any non-negative Borel function $\mu(t)$ (mutation rate at time $t$), we could define the mutation point measure as the Poisson point measure with intensity $\mu(d(\rho,x)) \, \lambda(dx)$. 
But if we are only interested in the genetic composition of the population of individuals/species co-existing at the same time, we can restrict our attention to ultrametric trees and in virtue of Theorem \ref{thrm:second sense}, focus on the tree $\tau_f(T)$ associated with a comb $f$ defined on the interval $I=[0, a]$. Let $\mu$ denote a diffuse Radon measure on $(0,T]$. From now on, we define the \textbf{mutation point measure} $M$ on $\text{Sk}$ as the Poisson point measure with intensity measure 
$$
\sum_{t\in[0,a]:f(t)\not=0}\delta_t(dx)\otimes \mu(dy)\indic{y<f(x)}, 
$$
where we assume that $f(0)=T$ so that there can also be mutations on the \textbf{origin branch} $L_0$. 
We will refer to the measure $\mu$ as the \textbf{mutation rate}. Usually, $\mu(dx)$ is taken equal to $\theta\,dx$, except when the comb studied is the image of an infinite tree by some map $\varphi$, as in Subsection \ref{subsec:yule}.

\subsection{Kingman coalescent with mutations}

Here, we consider that the comb $f$ is the Kingman comb and that $\mu(dx) =\theta\, dx$. The well-known results reviewed here are wonderfully exposed in \cite{DurBook2}.\\

Let us focus  first on the allelic partition of a sample of $n$ individuals as defined in the previous subsection. Let $A_n(k)$ denote the number of blocks of the allelic partition containing $k$ elements. Observe that we must have $\sum_{k=1}^n kA_n(k)=n$.
\begin{thrm}[Ewens sampling formula \cite{Ewe72}]
\label{thm:esf}
The random vector $(A_n(1),\ldots, A_n(n))$ has the same law as the random vector $(Y_1,\ldots, Y_n)$ conditional on $\sum_{k=1}^n kY_k=n$, where the $Y_k$'s are independent, and $Y_k$ is a Poisson r.v. with parameter $\theta/k$. For any vector $(a_1,\ldots, a_n)$ such that $\sum_{k=1}^n ka_k=n$,
$$
\xP(A_n(1)=a_1,\ldots, A_n(n) = a_n) = c_{\theta,n}\prod_{k=1}^n 
\frac{\left(\frac{\theta}{k}\right)^{a_k}}{a_k!} 
$$
where
$$
c_{\theta,n} := \frac{n!}{\theta(\theta+1)\cdots(\theta+n-1)} 
$$
\end{thrm}
\begin{figure}[!ht]
\includegraphics[width=10cm]{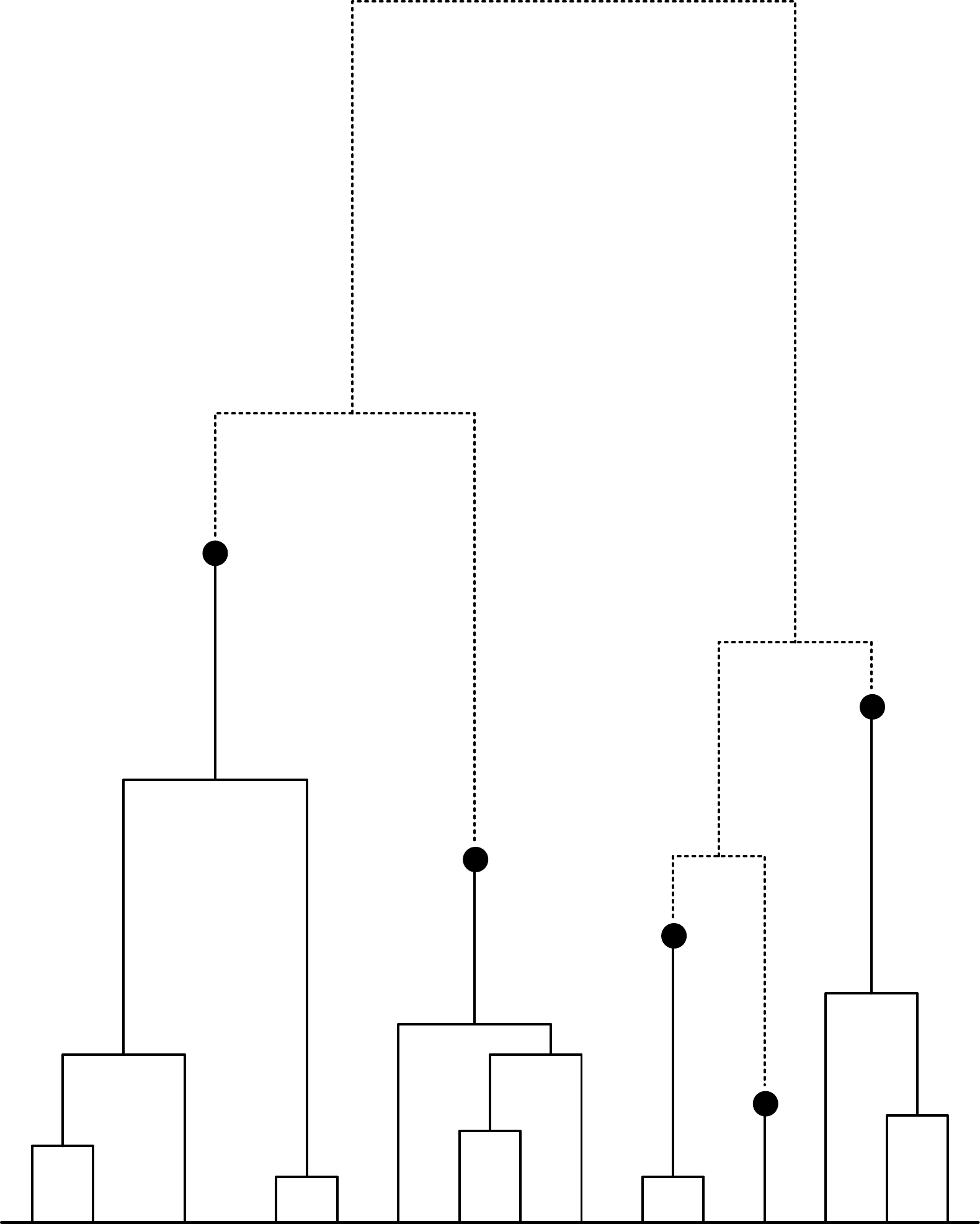}
\caption{Coalescent with mutations. Backward time flows upwards, black dots symbolize mutation events. Lineages further in the past of mutations are dotted to recall they are independent of the frequency spectrum, see proof of Theorem \ref{thm:esf}. In the present example, the resulting allele frequency spectrum is $(1, 1, 1, 1, 1, 0, 0...)$ (one block of each cardinality up to 5). }
\label{fig:esf}
\end{figure}
\begin{proof}
The trick is to follow the lineages of the $n$ individuals backwards in time. It is known (but not obvious from Definition \ref{dfn:kingman-comb}) that regardless of mutations, each pair of lineages coalesces at unit rate independently. Now each lineage is hit independently by a mutation at constant rate $\theta$. Once it is hit, mutations occurring further in the past can be ignored because they have no consequence on the allelic partition at present time, thanks to the infinite-allele assumption. In addition, by the sampling consistency of the Kingman coalescent, the lineage itself can be frozen and put aside because its presence has no consequence on the process of coalescences of other lineages. So we are left with a Markov process in backward time where pairs of lineages coalesce at unit rate and single lineages are frozen at rate $\theta$. Once all lineages have been frozen by a mutation, each of these mutations corresponds to a distinct allele and the maximum number of lineages downstream from this mutation is the size of the block corresponding to this allele. See Figure \ref{fig:esf} for an illustration.

Now it is not difficult to see that changing the arrow of time once again and changing time, this process becomes a pure-birth process with immigration started at 0 and stopped when it hits $n$.  It has unit birth rate and immigration rate $\theta$. So the allele frequency spectrum has the same law as the vector of sizes of immigrant families after the process is stopped.  For the calculations ending the proof see \cite{Lam11} or the references therein.
\end{proof}

\begin{crllr} 
\label{cor:hfs}
As $n\to\infty$,
the random vector  $(A_n(1),\ldots, A_n(n))$ converges in the sense of finite-dimensional distributions to a sequence $(Y_k)_{k\ge 1}$ of independent r.v., where $Y_k$ is a Poisson r.v. with parameter $\theta/k$.
\end{crllr}
\noindent
This limiting spectrum of \textbf{small families} (i.e., blocks with size $O(1)$) is sometimes called the \emph{harmonic frequency spectrum}. Besides the dependence in $1/k$, it is important to remember that the number of small families is $O(1)$, which is in deep contrast with what happens for coalescent point processes (see next subsection).\\

Now to understand the behavior of the \textbf{large families} (the other end of the spectrum), we set $X_n(i)$ the size of the $i$-th largest block of the allelic partition. The following result follows from the description of the allelic partition in terms of a pure-birth process with immigration, as in the proof of Theorem \ref{thm:esf}.
\begin{prpstn}[\cite{DT86}] 
The random vector $\left(n^{-1}X_n(k)\right)_{k\ge 1}$ converges in the sense of finite-dimensional distributions to the sequence $(P_k)_{k\ge 1}$ defined as
$$
P_k=Z_k\prod_{i=1}^{k-1}(1-Z_i),
$$
where the $(Z_i)$ are i.i.d. with density $\theta (1-z)^{\theta -1}$ (Beta $(1,\theta)$).
\end{prpstn}
\noindent
In words, the sizes of the largest blocks are of the order of the sample size $n$, and they represent fractions of the sample size which converge to the GEM distribution with parameter $\theta$, distribution that can be obtained by a recursive stick-breaking procedure of the unit interval. 

One could also tackle directly the problem on the Kingman comb, by considering the measures, say $X_k:=\ell (R_k)$, of the largest blocks $(R_k)_{k\ge 1}$ of the allelic partition of the whole population. One would theoretically expect $(X_k)$ to follow the GEM distribution. Another interesting question is to study the geometry of the subsets $(R_k)$ of the unit interval analogously to what is done below in Theorem \ref{thm_cpp_clonal}.

\subsection{CPP with mutations}

CPPs are less prone than the Kingman comb to calculations based on a sample of a fixed size. On the other hand, it can be shown that the genealogy of a sample of points thrown on the boundary according to a Poisson point process with constant intensity is equal in distribution to a new CPP, with finite intensity measure \cite{Lam09}. Modulo a modification of $\nu$, we can therefore focus on the allelic partition of the whole population without loss of generality. \\

From now on, we assume that the comb $f$ is a CPP killed at its first atom with second component larger than $T$, and with intensity $\nu$, a diffuse measure on $(0,\infty)$ such that  $\overline{\nu}(x):=\nu([x,\infty))<\infty$ for all $x>0$. We additionally assume that the mutation intensity measure $\mu$ is a diffuse Radon measure on $[0,\infty)$. In particular, we can define $\underline\mu(t):=\mu([0,t])$. Then it can be seen \cite{DL17} that if  $\int_0 \underline{\mu}(x)\nu(d x) < \infty$ the total number of mutations is finite a.s. whereas if $\int_0 \underline{\mu}(x)\nu(d x) = \infty$ then the number of mutations in any clade (set of all descendants of a point) is infinite a.s.\\

\noindent
Recall that a point of $I$ is clonal if it carries the same allele as the root, that is, there is no mutation on its lineage.
The next statement is based on the idea that for any clonal $t$, the right-hand side of $t$ only `sees' the mutation-free lineage of $t$.   
\begin{thrm}[\cite{DL17}] 
\label{thm_cpp_clonal} 
Let $T=\infty$ so that $I=[0,\infty)$ and let $R$ denote the closure of the clonal set. Conditional on the absence of mutation on the origin branch $L_0$, $R$ is a regenerative set that can be described as the range of a subordinator whose Laplace exponent $\phi$ is given by
        \[\frac{1}{\phi(\lambda)} = \int_{(0, \infty)} \frac{e^{-\underline{\mu}(x)}}{\lambda + \overline{\nu}(x)} \mu(d x).\]
        \end{thrm}
        \noindent
In addition, there exist semi-explicit formulae (that we will not provide explicitly here) for the \textbf{expectation of the allele frequency spectrum}, namely for
$$
\Lambda_T(dx):= \E(A_T(dx)),
$$
where $T$ is here to recall the dependence on the height of the CPP, and the measure $\Lambda_T$ is defined on $\xN$ when $\nu$ is finite \cite{Lam09} and on $(0,\infty)$ when $\nu$ is infinite \cite{DL17}. In addition, the following convergence holds weakly on $[0,\infty)$ 
$$
\Lambda(dx):=\lim_{T\to\infty}\frac{\Lambda_T(dx)}{\E(a(T))},
$$
where $a(T)$ is the width of the CPP killed at its first atom with second component larger than $T$ (in particular, $\E(a(T))=1/\nu([T,\infty))$). In the case of a unit rate critical birth-death process with mutations at rate $\theta$, we get \cite{Lam09} 
$$
\Lambda(\{k\}) = \frac\theta k \,(1+\theta)^{-k}\qquad k\ge 1.
$$
Analogously, in the case of the Brownian CPP with intensity measure $\frac{dx}{x^2}$ and mutation measure $\theta\,dx$, elementary calculations show that
$$
\Lambda (dx) =\frac{\theta}{x}\,e^{-\theta x} dx\qquad x>0.
$$
The last two formulae are reminiscent of the harmonic frequency spectrum displayed in Corollary \ref{cor:hfs}.\\

Next, a LLN type of argument entails the following result.
\begin{thrm}[\cite{CL12, DL17}]
\label{thm:LLN}
Let $\bar A_T(q)$ denote the number of alleles whose carriers on the boundary form a set of measure larger than $q$ (counting measure in the finite case, $\ell$ otherwise). Then the following limit holds a.s.
$$
\lim_{T\to\infty}\frac{\bar A_T(q)}{a(T)} = \Lambda([q,\infty)).
$$
\end{thrm}
\noindent
A similar convergence will be mentioned in the next subsection for the allele frequency spectrum at time $t$ of the supercritical branching process  as $t\to\infty$.

\subsection{Supercritical processes with mutations}

Here, we wish to explain how the allelic partition at the boundary of a supercritical branching process can be understood from the allelic partition of a CPP, in the vein of Theorem \ref{thm_yule}. We also say a word on the limiting allele frequency spectrum for large times. \\

Let $\tr$ be the infinite tree generated by a supercritical birth-death process with birth intensity measure $\beta$ conditioned on non extinction and let $q(t)$ denote the probability of extinction of a particle born at time $t$ (which depends on some unspecified, possibly 0, possibly age-dependent death rate). Also assume that conditional on $\tr$, mutations occur on the lineages of $\tr$ according to some intensity measure $\mu$ on $[0,\infty)$. Also recall that there is a measure $\mathscr L$ defined on the boundary of $\tr$ by 
 \[
   \mathscr{L}\left (B_u\right ) := \lim _{t\uparrow \infty} \frac{\tilde N_u(t)}{e^{\underline{\tilde\beta}(t)}},
  \]
where $B_u$ denotes the set of infinite sequences $v$ with prefix $u$, $\tilde N_u(t)$ is the number of descendants of $u$ at time $t$ which have infinite descendance and $\underline{\tilde\beta}(t) = \int_0^t (1-q(s))\,\beta(ds)$.

Now from Theorem \ref{thm_yule} and the remarks following it, if $\varphi:[0,\infty)\to (0,T]$ is a decreasing bijection, then $\varphi(\tr)$ is isometric to the tree $\tau_f$ constructed from a CPP $f$ with intensity $\nu$ and height $T$, naturally measured by $\ell$, where
$$
\nu([t,\infty)) = \exp\left(\underline{\tilde\beta}\circ\varphi^{-1}(t) \right)\qquad t\in[0,T].
$$
Now it is straightforward to see that the mutations on $\tau_f$ occur according to a mutation point measure with mutation rate $\mu_\varphi:=\mu\circ\varphi^{-1}$.\\ 

If we assume that \textbf{$\mu$ has finite mass $\theta$} and we take $\varphi$ given by $\varphi(t) =\theta^{-1}\mu([t,\infty))$, then the allelic partition at the boundary is equal in distribution to that of a CPP (with intensity $\nu$ previously displayed) with height 1 and mutations \textbf{at constant rate $\theta$}.\\

On the contrary, if we take $\mu(dx) = \theta \,dx$ and $\varphi$ given by $\varphi(t) =e^{-\underline{\tilde\beta}(t)}$, then $\varphi(\tr)$ has the distribution of $\tau_f$ where $f$ is a Brownian CPP with height 1, intensity $\frac{d x}{x^2}$ and mutation rate $\mu_\varphi$ satisfying $\mu_\varphi[(x,1]) = \theta \varphi^{-1}(x)$.
In particular, if the birth-death process is time-homogeneous with birth rate $b$ and if we set $a:= b(1-q)$, then $\varphi(t) = e^{-at}$, so that $\varphi^{-1}(t) =-\ln(t)/a$. As a conclusion, the allelic partition at the boundary of a time-homogeneous supercritical birth-death process with mutations \textbf{at constant rate $\theta$} is equal in distribution to that of a \textbf{Brownian CPP} with height 1, intensity $\frac{d x}{x^2}$ and mutations \textbf{at inhomogeneous rate $ \theta/(at)$}.\\

\noindent Similarly to Theorem \ref{thm:LLN}, the allele frequency spectrum at time $t$ of the supercritical birth-death process properly rescaled converges a.s. as $t\to\infty$. The LLN type of argument invoked here is known as the theory of branching processes \textbf{counted with random characteristics} \cite{Ner81, JN84a, JN84b, JN96}. In our setting, the random characteristic of individual $i$, say, can be for example the number $\chi_i^{k}(t)$ of mutations that $i$ has experienced during her lifetime and which are carried by  $k$ alive individuals, $t$ units of time after her birth ($\chi_i(t)=0$ if $t<0$). Then the total number $A_t(k)$ of alleles carried by $k$ individuals at time $t$ (except possibly the ancestral type) is the sum $\sum_i \chi_i(t-\alpha_i)$ over all individuals $i$ (dead or alive), where $\alpha_i$ is the birth time of individual $i$. The theory of branching processes counted with random characteristics ensures that these sums rescaled by $e^{at}$ converge a.s.\ on the survival event.

This method has been used extensively in \cite{TaiBook}. Further refinements have been obtained in \cite{CL13} thanks to the use of coalescent point processes. For example, we have shown the convergence of the (properly rescaled) largest blocks of the partition at time $t$, as $t\to\infty$. In particular, we showed, letting $N\equiv N_t$ denote the total population size at time $t$, that the largest blocks are roughly of size $N^{\theta/a}$ when $\theta<a$, $\log(N)^2$ when $a=\theta$ and $\log(N)$ when $\theta >a$. Similarly, we showed that the oldest mutations with alive carriers at time $t$ appeared roughly at time $O(1)$ when $\theta <a$, at time $\ln(t)/a$ when $\theta =a$ and at time $(1-a/\theta)t$ when $a<\theta$. Notice that in contrast to the Kingman case, no block has a size of the order of the total population size when the mutation rate is constant. It is an open question to investigate whether the compactification mentioned earlier in this subsection of the tree into a CPP can shed extra light on these results. In particular, the case of large families is not straightforward at all since the compactification focuses precisely on subsets of the boundary with measure of the order of the total population size.

\bibliographystyle{plain}
\bibliography{RefsMAS}

\end{document}